\documentclass[12pt]{amsart}
\usepackage{amsfonts,amsthm,amsmath,amssymb}
\usepackage{comment}

\def\IH{{{\mathbb{H}}^3}}

\def\cal{\mathcal}

\def\Hex{\cal{H}}
\def\Pent{\cal{P}}

\def\oC{{\hat{\mathbb {C}}}}

\def\tr{{\rm tr}}

\def\IR{{\Bbb R}}
\def\IP{{\Bbb P}}

\def\IT{{\Bbb T}}

\def\IC{\Bbb C}

\def\x(AB)i{{Ax_{{A \cdot B^{-1}}}}}
\def\(AB)i{{{{A \cdot B^{-1}}}}}

\newtheorem{thm}{Theorem}[section]

\newtheorem{lemma}[thm]{Lemma}

\newtheorem{cor}[thm]{Corollary}
\newtheorem{definition}[thm]{Definition}
\newtheorem{prop}[thm]{Proposition}

\begin{document}
\title{Planar families of discrete groups}
\author{ Jane Gilman}
\address{ Department of Mathematics, Smith Hall\\
 Rutgers University, Newark, NJ 07102 \\ }
\email{gilman@andromeda.rutgers.edu}
\thanks{Supported in part by NSA grant MSPF\#02G-186.}
\author{ Linda Keen}
\address{ Mathematics Department\\CUNY
Lehman College and Graduate Center \\ Bronx, NY 10468, USA}
\email{linda.keen@lehman.cuny.edu}
\thanks{ Supported in
part by  PSC-CUNY grant}

\subjclass{Primary 30F10,30F35,30F40; Secondary 14H30, 22E40}
\keywords{Kleinian group, Schottky group, M\"obius group, discrete
group}

\maketitle
\section{Introduction}
Determining the space of free discrete two generator groups of
M\"obius transformations is an old and difficult problem. In this
paper we show how to construct large balls of full dimension in
this space.
 To do
this, we begin with a marked discrete group of non-separating
disjoint circle type, {\em nsdc group} (see section~\ref{sec:nsdc}
for the definition). Such a group determines three disjoint or
tangent planes. We prove that there is a whole family of discrete
groups, which we term the {\em nsdc-planar family}, that share
these planes.

From  specific information about the nsdc planes which depends on
the marking data, we show how
 to construct all the discrete groups in the nsdc planar family.
In particular, we find a set of  six real numbers that serve as
parameters  for this family.  We then construct an embedding of
our parameters into a classical representation of the full space
of free discrete groups
as a subset of $\oC^3$.  We see that each planar family fills out
a ball of full dimension in the classical embedding.

We remark that our construction does not use the usual
theory of quasi-conformal deformations of a given group nor does
it depend simply on the use of coverings and quotients that give
commensurable groups.

The paper is organized as follows. In section \ref{sec:not}  we
give the basic notation and definitions and in section
\ref{sec:planar}  we define the nsdc planar family. Sections
\ref{sec:ortho} and \ref{sec:perp} develop computational
techniques and in section \ref{sec:constr} we construct our new
parameters. In section \ref{sec:moduli} we show how to our
parameters are related to the classical parameters. Finally, in
section \ref{sec:neccs} we give a necessary and sufficient
condition for a marked group to belong to a given planar family.
\section{Notation and definitions} \label{sec:not}
An element $X$ of $PSL(2, \IC)$ acts as a M\"obius transformation
on $\oC$, the complex sphere. The  action extends in a natural way
to hyperbolic three-space $\IH$. When considered as the boundary
of hyperbolic three-space, $\oC$ is referred to as the sphere  at
infinity.

Elements of $PSL(2, \IC)$ are classified as
  loxodromic, elliptic, parabolic according to the square of their
traces. The classification of transformations can also be
described by the  action  on $ \oC$ or on $ \oC \cup \IH$.

We follow the notation of
\cite{Fench}.
For $x \mbox{ and }  y \in (\oC \cup \IH)$ with $x \ne y$ we
let $[x,y]$ denote the
oriented hyperbolic geodesic  passing through $x$ and $y$.
 The {\sl ends} of the
geodesic $[x,y] $ are by definition the points $v$ and $w$
with $\{ v, w \} =  [x,y] \cap \oC$.
If $v$ and $w$ are the ends of $[x,y]$, then $[x,y]=[v,w]$.
 The notation $[x,y]$  also indicates a direction  so that
$[y,x]$ is the geodesic with the opposite orientation.

Following Fenchel \cite{Fench}, we include {\sl improper lines}, in our
considerations. A {\sl proper line} in one whose ends are distinct.
 An improper line is one whose ends coincide.

 A M\"obius transformation fixes one or two points on $\oC$.
A loxodromic transformation $A$ has two distinct fixed points as does
an elliptic transformation. These transformations
fix the geodesic $X$
in $\IH$ whose ends are the fixed points of $A$ on $\oC$. This geodesic is called the {\sl axis} of $A$ and is denoted  by $Ax_A$.

A parabolic transformation has one fixed
point. Using the terminology of improper lines, parabolic transformations
also have axes. If $A$ is a parabolic
transformation with fixed point $u \in \oC$, we consider $[u,u]$
to be its axis.

Every proper pair of lines in $\IH$ has a unique common perpendicular.
If  $u,v,w,z$ are points in $\oC$, then two distinct improper lines
$[u,u]$ and $[v,v]$ have the unique common perpendicular $[u,v]$.
The common perpendicular to the  improper line $[u,u]$ and the proper
line $[v,w]$ is  the geodesic
through $u$ perpendicular to $[v,w]$.

With this terminology,  every pair of elements $A, B \in
 PSL(2,\IC)$ determines a unique hyperbolic line in $\IH$, the
 common perpendicular of their axes; we denote this common perpendicular
by $L$
 and denote the ends of $L$ by $n$ and $n'$.

\subsection{Half-turns}\label{sec:halfturns}
For any hyperbolic line $[x,y]$ we let $H_{[x,y]}$ be the
half-turn about the line with ends $x$ and $y$. We note
that $H_{[x,y]}$ fixes every hyperbolic plane $\IP$
whose horizon  $C_P$ passes through $x$ and $y$ (\cite{GiNSDC}).
The half-turn will
interchange the exterior and the interior of $\IP \in \IH$.
\subsection{NSDC Groups.} \label{sec:nsdc}
There is a natural way to attach  a
three generator group and six complex numbers to any a two generator group (see \cite{GiNSDC}). 

Let $G = \langle A,B \rangle$ and let  $L$ be
 the common perpendicular to  $A$ and $B$.
 There are unique hyperbolic lines $L_A$ and $L_B$ such that  $A =  H_{L_A}
\cdot H_L$ and $B =  H_{L_B} \cdot H_L$. We let $a$ and $a'$ be
the ends of $L_A$ so that  $L_A = [a,a']$ and let $b$ and $b'$ be those
of $L_B$ so that $L_B = [b,b']$.   We define
\begin{definition}
The marked {\sl three generator group} determined by
 $G = \langle A,B \rangle$ is denoted by $\IT G$ and defined by
$\IT G = \langle H_{L_A}, H_L, H_{L_B}
\rangle$.
\end{definition}
By construction $G$ is a normal subgroup of $\IT G$ of index at most two which immediately
implies
\begin{prop}\label{prop:index2} $G$ is discrete if and only if $\IT G$ is discrete. \end{prop}

We also define
\begin{definition}
The {\sl ortho-end} of $G$ is the six-tuple of complex numbers
$(a,a',n,n',b,b')$
\end{definition}
\begin{definition} Six points in $(a,a',n,n',b,b') \in \oC^6$ have the
{\sl non-separating disjoint circle property} if there are pairwise disjoint or  tangent
circles
on $\oC$,  $C_A$, $C_D$ and  $C_B$ (respectively)
 passing through $a$ and $a'$,
$n$ and $n'$, and $b$ and $b'$ (respectively) such that no one circle separates
the other two.
\end{definition}
\begin{definition} A marked group $G =\langle A, B \rangle$
 is a marked {\sl  non-separating disjoint circle group or nsdc group}
if the ortho-end of  $A$ and $B$ has the non-separating disjoint circle
property.
$G$ is a {\sl  non-separating disjoint circle group }
if some pair of generators for $G$ has the non-separating disjoint circle
property. The corresponding group $\IT G$ is also called an {\sl nsdc group}.
\end{definition}
It was shown  in \cite{GiNSDC} that
\begin{prop} \label{prop:discrete} Let $G$ be a two-generator subgroup of $PSL(2, \IC)$.
If some  ortho-end of $G$ has the non-separating disjoint circle property,
then $G$ is discrete.
\end{prop}
\section{The planar  family of  nsdc groups} \label{sec:planar}
Let $G$ be   an nsdc group with ortho-end $(a,a',n,n',b,b')$. Let
the three  circles $C_A$, $C_D$, $C_B$ be a set of non-separating
disjoint circles for this ortho-end. Let $\IP_A$, $\IP_D$ and
$\IP_B$ be planes in $\IH$ whose respective horizons are $C_A$,
$C_D$, and $C_B$. Let $L_A'$ be any hyperbolic line lying on
$\IP_A$, let $L_B'$ be any hyperbolic line lying on $ \IP_B$, and
let $L'$ be any hyperbolic line lying on $\IP_D$. Set $\IT G' =
\langle H_{L_A'}, H_{L'}, H_{L_B'} \rangle$ and $G'=\langle
H_{L_A'}\cdot H_{L'},H_{L_B'}\cdot H_{L'}\rangle$.
\begin{prop}
$\IT G'$ and $G'$ are both discrete.
\end{prop}
\begin{proof}
Let $[\alpha,\alpha']$,$[\eta,\eta'],[\beta,\beta']$ be the
ends of $L_{A'}$, $L'$ and $L_{B'}$ respectively.  The six-tuple
$(\alpha,\alpha',\eta,\eta',\beta,\beta')$ is an ortho-end of $G'$
and $\IT G'$ with nsdc circles $C_A,C_D$ and $C_B$.  By
proposition~\ref{prop:discrete} $G'$ is discrete and by
proposition~\ref{prop:index2} $\IT G'$ is also.
\end{proof}

The family of groups $\IT G'$ depends upon $\IP_A$, $\IP_D$ and
$\IP_B$. We write $\IP$ for $\IP_D$.

\begin{definition} We call the triple of planes
$(\IP_A, \IP, \IP_B)$ a {\sl non-separating disjoint planar
triple} or an {\sl nsdc}-triple if the horizons are a set of {\sl
nsdc} circles. We call the set of all nsdc groups with a fixed
planar triple,  a {\sl planar family of nsdc groups}.
\end{definition}

 We want to
explore all discrete groups corresponding to a given nsdc planar
triple.
\section{Ortho-ends and pull back circles} \label{sec:ortho}
In this section we present notation needed to describe
a triple of nsdc planes.
We begin with $G= \langle A, B \rangle$
and assume that the ortho-end of $G$ is the six-tuple
$(a,a',n,n',b,b')$ and that the nsdc-planes are $\IP_A, \IP, \IP_B$
respectively. The planes each are determined by their {\sl pull-back angles},
which can be described as follows:

Let  $k$ and $k'$ be  any two
points in $\oC$. Any circle through $k$ and $k'$
has center lying on the perpendicular bisector of the line connecting
$k$ and $k'$. Thus its  center is at
$c_t = {\frac{k+k'}{2}} + it{\frac{k-k'}{2}}$ for some real number $t$.
We think of the  center as having been pulled back $t$ units from
${\frac{k+k'}{2}}$.
We let $\theta_t$ denote the angle from
the radius connecting $c_t$ to $k$ to
the Euclidean segment joining $k$ and $k'$ so that
$t = {\frac{|k-k'|}{2}}  \cdot \tan \theta_t$.  We call
$\theta_t$ the {\sl pull-back angle} of the circle.
Each circle passing through $k$ and $k'$ corresponds to a
unique pull-back angle $\theta_t$ with  $-\pi/2 \le \theta_t <  \pi/2$.
If $r(t)$ denotes the radius of the circle with pull back angle
$\theta_t$, we have $r(t) = {\frac{{|k-k'|} }{2\cos
\theta_t}}$.

When we are discussing a fixed plane whose horizon passes through
$k$ and $k'$, we write $\theta_K$
for the pull-back angle. We write $r_{\theta_K}$ for
the radius of the horizon of the plane  and $c_{\theta_K}$ for its center.

In particular, if we have a family of nsdc-planes for
$G = \langle A,B \rangle$,
we assume that
 $\IP_A$ has pull-back angle $\theta_A$,
 $\IP_B$ has pull-back angle $\theta_B$,
and  $\IP$ has pull-back angle $\theta$.

\section{Computations and perpendiculars} \label{sec:perp} In
what follows we will be interested in being able to compute
quantities such as points and planes directly from
 the  ortho-end $(a,a',n,n',$ $b,b')$ and the three pull-back
angles ($\theta_A, \theta, \theta_B$).

We define three points $v \in \IP$, $v_A \in \IP_A$ and $v_B \in
\IP_B$ by $v = L \cap Ax_A$,  $v_A = L_A \cap Ax_A$, $v_B = L_B
\cap Ax_B$. These are  all computable from the ortho-end and the pull-back angles.

If $\IP$ is any hyperbolic plane,  $L^{\IP}$ a hyperbolic line on
$\IP$ and $x$ a point on $L^{\IP}$, there is a line perpendicular
to $\IP$ passing through $x$ which we denote by $V^{\IP}_x$ and a
unique line perpendicular to $L^{\IP}$ and $V^{\IP}_x$ lying on
$\IP$ and passing through $x$ which we denote $M^{\IP}_x$. If
there is no confusion we will omit the superscript $\IP$ from the
notation for these lines.
Let $L_W$ be a line with ends $[w,w']$ and suppose it lies on a
plane $\IP_W$ with pull-back angle $\theta_W$. Then we can compute
that $$c_{\theta_W}= {\frac{w+w'}{2}} + i{\frac{(w-w')^2}{4}}\tan
\theta_W.$$ Thus for any $x$ on $L_W$
$$P^{\IP_W}_{x} = [ {\frac{w+w'}{2}} + i{\frac{(w-w')^2}{4}}\tan \theta_W,
x].$$ Also since one  end of $M^{\IP_W}_x$ must be
${\frac{w+w'}{2}} + i(r_{\theta_W})$, we have
$$M^{\IP_W}_x = [{\frac{w+w'}{2}} + i(r_{\theta_W}), x].$$
The point is that this
explicit calculation only depends upon $\theta_W$. We note that a
line perpendicular to $\IP_W$ passing through a point $x$ in
$\IP_W$ is the line $[c_W, x]$. We are interested in these
perpendiculars when $x= v_A$.
\section{Constructing the planar family}\label{sec:constr}
A loxodromic transformation that is a pure translation by distance $d$ with zero rotation angle  is called a
{\em hyperbolic transformation with translation length $d$}.
\begin{definition} Let $X$ be any hyperbolic line and $\tau$ any angle.
We let $R_{X,\tau}$ be the elliptic transformation that is  rotation through an angle $\tau$
about the
line $X$. If $d$ is any positive real number, we let $T_{X,d}$
be the hyperbolic transformation whose  axis is $X$ and whose
translation length is $d$.
\end{definition}
We  let our rotation angles  lie in the interval $(-\pi,\pi]$ in
order to keep track of the orientation of a line under the
rotation.

We note that given any nsdc planar family, these two types of
moves, rotation about a line perpendicular to a plane of the
family and parallel transport along a line lying in a  plane of
the family preserve the family. Specifically we prove

\begin{lemma} Assume that $G = \langle A, B \rangle$ is of nsdc type in the $(\IP_A, \IP, \IP_B)$   planar family.
Let $L_{A'}= g(L_A)$ where
  $g= g_1,...,g_n$ and  for each $i$, either
\begin{itemize}
\item  $g_i= T_{X_i,d_i}$ where  $X_i \in \IP_A$ and $d_i$ is a real number
or
\item $g_i = R_{X_i, \tau_i}$ where $X_i$ is perpendicular to $\IP_A$ and
$\tau_i$ is any angle $-\pi < \tau_i \le \pi$.
\end{itemize}
Then
$\langle H_{L_{A'}}, H_L, H_{L_B} \rangle$ is a discrete group
in the $(\IP_A, \IP, \IP_B)$   planar family.
If  $A' = H_{L_A'}H_L$, then
$G'= \langle A', B \rangle$ is a discrete group of nsdc type.
\end{lemma}
\begin{proof} By construction each of the $g_i$ maps lines on $\IP_A$ to
lines on $\IP_A$; $g_i H_{L_A} g_i^{-1} = H_{g_i(L_A)}$.
\end{proof}
Thus these moves preserve the planar family. Next we see that
every triple of geodesics $(L_{A'},L', L_{B'})$ lying in the
respective planes   $\IP_A, \IP, \IP_B$  of the planar family can
be so obtained. Precisely,
\begin{prop} \label{prop:dtau}
Assume that  $G = \langle A, B \rangle $ is of nsdc type with
ortho-end $(a,a',n,n',b,b')$ and  planar family $(\IP_A, \IP, \IP_B)$ with
pull back-angles $(\theta_A, \theta, \theta_B)$. Let $v_A = L_A \cap
Ax_A$. Then   $G'= \langle A', B \rangle $ is in the nsdc planar
family with $A' = H_{L_{A'}}\cdot H_L$ where  $L_{A'}$ is obtained from
$L_A$ by the following (where the plane $\IP_A$ is implied).
\begin{enumerate}
\item $X =
R_{V_{v_A}, \tau }(M_{v_A})$ where $\tau$ is an  angle
with $-\pi \le \tau < \pi$,
\item $Y = R_{V_{v_A}, \tau }(L_A)$, and
\item $L_{A'} = T_{X, d}(Y)$ for some positive number $d$.
\end{enumerate}
Conversely,  given any $L_{A'} \in \IP_A$ there exists an angle
$\tau$,   $-\pi \le \tau < \pi$ and a positive number $d$ so that
$L_{A'}$ can be obtained from $L_A$ by the above construction.
\end{prop}
\begin{proof}
If $L_{A'}$ is any line on $\IP_A$, we can drop a perpendicular
from $v_A$ to $L_{A'}$, whether or not the lines $L_A$ and $L_{A'}$
are disjoint. Let $X$ be the geodesic on which this perpendicular
lies, oriented towards $v_A$, and let $d$ be the distance from  $L_{A'}$ to $v_A$ along $X$.
Then $T^{-1}_{X,d}(L_{A'}) = Y$ passes through $v_A$ and makes
some angle $\tau$ with $L_A$. We can rotate  $Y$ to  $L_A$ by
applying $R_{V_{v_A},\tau^{-1}}$. (Alternatively, rotating by
$\pi/2+\tau$  sends  $Y$ to $M_{v_A}$.)  Combining these moves we
have
$$L_{A'}= (T_{R_{V_{v_A}, \tau }(L_A), d} \circ R_{V_{v_A},\tau})(L_A)
= (T_{R_{V_{v_A}, \pi/2 + \tau }(M_{v_A}), d} \circ R_{V_{v_A},\tau})(L_A)$$

Thus any line $L_A'$ on $\IP_A$  determines a positive number
$d=d_A$ and an angle $\tau=\tau_A$ in $(-\pi,\pi)$. Conversely, given the
initial data  $(\theta_A, a \mbox{ and } a')$ and the pair $(d_A,
\tau_A)$ we can use this construction to find $L_{A'}$.
\end{proof}
\begin{definition} For the fixed planar family $(\IP_A, \IP, \IP_B)$ with pull back-angles $(\theta_A, \theta, \theta_B)$ we call the sequence of moves in the proposition {\sl the $(d_A,\tau_A)$
moves} determining $L_A'$. \end{definition}
\begin{cor}
All of these moves can be calculated directly from
$(a,a',n,n',b,b)$ and $(\theta_A, \theta, \theta_B)$.
\end{cor}
\begin{proof}
This is the content of section \ref{sec:perp}  on perpendiculars.
\end{proof}

In summary we have
\begin{thm} Let $G = \langle A,  B\rangle$ be a fixed marked group of nsdc type in the planar family $(\IP_A, \IP, \IP_B)$ with pull back-angles $(\theta_A, \theta, \theta_B)$.  Then the set of triples
$${\mathcal D}=\{(d_A,\tau_A),(d,\tau),
(d_B,\tau_B) \in ([0,\infty),(-\pi,\pi])^3\}$$ form a set of
moduli for the family. That is, given an ortho-end
$(a,a',n,n',b,b')$ in $\oC^6$ and pull-back angles $(\theta_A,
\theta, \theta_B)$ there are six real numbers $(d_A,\tau_A$, $d,
\tau$,$ d_B,  \tau_B)$  such that every nsdc group $G' = \langle
A', B' \rangle$ in the planar family $(\IP_A, \IP, \IP_B)$ can be
obtained through the moves $(d_A,\tau_A)$, $(d,\tau)$  and
$(d_B,\tau_B)$ giving the lines $L_{A'}$,   $L'$ and $L_{B'}$
respectively and hence the marked group $G'$.
\end{thm}

\begin{proof} Proposition \ref{prop:dtau} shows that given each pair $(d_A,\tau_A)$,
appropriate moves can be chosen to find $L_{A'}$ and conversely
each $L_{A'}$ determines such a pair; the same is true for $L'$
and $L_{B'}$.  Thus all marked groups with these nsdc planes  are
determined. The point is that given the planar triple, or equivalently
the three angles $(\theta_A, \theta_B, \theta)$, every
ortho-end that gives a group in this  planar nsdc family can be
obtained.
\end{proof}

\section{Relation to classical moduli}\label{sec:moduli}

The marked nsdc groups are a subset of the space ${\mathcal S}$ of marked discrete free groups $G=\langle A,B \rangle$.  Classically,  ${\mathcal S}$ can be embedded into a subset of $\IC^3$ by choosing as moduli $(\tr A, \tr B, \tr AB^{-1})$. Specifically, one picks a marked group $G$ as base point and arbitrarily chooses elements in $SL(2,\oC)$ (again called $A$ and $B$) to represent the generators.  The signs of the traces of these elements depends on this choice but once it is made, there is a unique matrix corresponding to every  other group element.  The space ${\mathcal S}$ is then obtained by quasiconformal deformation. The computation of the full boundary of ${\mathcal S}$ in this embedding is an open hard question.

 In this section we show how these trace parameters are related to the nsdc moduli we found above.

\subsection{Skew Hexagons}
Let $G = \langle A,  B\rangle$ be a fixed group of nsdc type in
the planar family $(\IP_A, \IP, \IP_B)$ with pull back-angles
$(\theta_A, \theta, \theta_B)$ and let $\IT G$ be the corresponding
group.  As we saw in section~\ref{sec:nsdc}, we can write $A =
H_{L_A} \cdot H_L$,  $B =  H_{L_B} \cdot H_L$ and
$AB^{-1}=H_{L_A}\cdot H_{L_B}$. Here we assume the axes are
oriented from their repelling fixed points to their attracting
fixed points. We assume $L$ is oriented from $Ax_A$ to $Ax_B$,
$L_A$ is oriented from $Ax_A$ to $Ax_{AB^{-1}}$ and $L_B$ is
oriented from $Ax_{AB^{-1}}$ to $Ax_B$.

Following Fenchel, \cite{Fench}, we can form a {\em skew right angled hexagon} associated to $G$, which we
denote by $\Hex_G$ with sides:

 $$\Hex_G\mbox{  has sides } L_A,Ax_A, L, Ax_B,L_B,Ax_{AB^{-1}}$$

 Here we adopt the convention that we label
sides of a hexagon by the hyperbolic line upon which they lie
taking the segment indicated by the order. That is,  six ordered
geodesics determine six vertices: if the geodesics are
$\tilde{s}_1, \tilde{s}_2, \tilde{s}_3, \tilde{s}_4, \tilde{s}_5,
\tilde{s}_6$,  set  $v_i= \tilde{s}_i \cap \tilde{s}_{i+1}$ for $i
\le i \le 6$ with indices taken modulo 6 and let the hexagon side
denoted by $s_i$ be the segment of $\tilde{s}_i$ traversed from
$v_{i-1}$ to $v_i$ where $v_0=v_6$.
With this convention it is
clear which segment of $\tilde{s}$ we are talking about, and,
therefore, we do not distinguish notationally between the segment
$s$ and the geodesic ${\tilde{s}}$ on which it lies. We write $s$
for both.

The labeling gives an orientation to the hexagon and its sides.
Note that the orientation of a side within the hexagon may be
opposite to the orientation of the geodesic containing the side.

 Relative to the  orientation of the hexagon
we define the {\em complex length} $\delta_i=\delta(s_i)$ {\em of
the side $s_i$} as follows: $t_i=\Re \delta(s_i)$ is the
hyperbolic length of $s_i$; $\psi_i=\Im \delta(s_i)$ is the angle
from $T_{s_i,t_i}(s_{i-1})$ to $s_{i+1}$ with indices taken modulo
6. Here $\Re$ and $\Im$ denote the real and the imaginary parts of
a complex number. The complex lengths satisfy the ``cosine rule''
$$ \cosh(\delta_{i+4})=\cosh(\delta_i)\cosh(\delta_{i+2})+
\cosh(\delta_{i+1})\sinh(\delta_i)\sinh(\delta_{i+2}).$$

It follows that the complex lengths of any three
alternating sides or any three adjacent sides determine the other
three lengths.

For example, in the hexagon $\Hex_A$ we have
$\delta(Ax_A)=t_A + i \psi_A$ where
 $t_A$ is the hyperbolic distance along $Ax_A$ between $L_A$ and  $L$ and $\psi_A$ is the angle from  $T_{Ax_A,t_A}(L_A)$ to $L$. ( Define $t_B,\psi_B$ and $t_{AB^{-1}}=t,\psi_{AB^{-1}}=\psi$ similarly.)
 It is a standard computation that $\tr A = 2 \cosh (t_A + i \psi_A)$ and similarly for $\tr B$ and $\tr AB^{-1}$.

Using the cosine rule, these three traces determine the complex lengths of the other three sides, $L_A,L,L_B$, of the hexagon and the complex lengths of the $L_A,L,L_B$ sides determine the traces.

\subsection{Hexagons for the planar family}

Let $G'=\langle A',B\rangle$ be a group in the planar family  where $L_{A'}$ is obtained from $L_A$ by $(d_A,\tau_A)$ moves.
We can form the pentagon $\Pent_{A,A'}$ with ordered sides
$Ax_A,L,Ax_{A'},$ $L_{A'}, X_{v_A} $ where $X_{v_A}$ is the
perpendicular from $v_A$ to $L_{A'}$.  We may consider the pentagon as a degenerate skew hexagon where the length of the (degenerate) fifth side has no real part as follows:
 $$s_1 = Ax_A; \quad \delta_1=t_A+i\psi_A$$
 $$s_2 = L; \quad \delta_2=t_L + \psi_L$$
 $$s_3 = -Ax_{A'}; \quad \delta_3= t_{A'}+i(\psi_{A'}+\pi)$$
 $$s_4=L_{A'}; \quad  \delta_4 = t_{L_A'}+ i \psi_{L_A'} $$
 $$s_5;  \quad \quad \delta_5= -i\tau_A$$
 $$s_6=X_{v_A}; \quad \delta_6=d_A$$
\begin{thm} Given the fixed marked group $G$ and the parameters
$(d_A,\tau_A)$, $(d,\tau)$, $(d_B,\tau_B)$ we can find the traces
of the elements $A',B',$ $A'B'^{-1}$ of  any group  in the planar
family $(\IP_A,\IP,\IP_B)$.
\end{thm}

\begin{proof} We use the cosine rule applied to the
degenerate skew hexagon above to find $\delta_3$ and hence the trace of $A'$. We get
$A'$ from its axis and its trace  and thus obtain the
group $G'=\langle A',B \rangle$.  We can also find the complex
distance $\delta_2$ along the  side $L$.  Since we know
complex distance along $L$ from $Ax_A$ to $Ax_B$ in $\Hex_G$, we can compute the
distance along $L$ from $Ax_A'$ to $Ax_B$ in $\Hex_{G'}$. Using the cosine rule again, we get all
the complex lengths in the hexagon $\Hex_{G'}$ and hence the traces $(\tr A', \tr B, \tr A'B)$.

 Since we can get to any group in the family by a sequence of moves as above, we can find its trace moduli as well.
\end{proof}

As an immediate corollary we have

\begin{cor} The constructions above determine an embedding of ${\mathcal C}$ into ${\mathcal S}$. Thus each planar family fills out a ball in ${\mathcal S}$ of full dimension.
\end{cor}

\section{Ortho-ends in the planar family }\label{sec:neccs}
Given a planar nsdc family with ortho-end $(a,a',n,n',b,b')$ and
planes $(\IP_a,\IP,\IP_B)$ with pull back angles $\theta_A,
\theta, \theta_B$, we can find necessary and sufficient conditions
for a six-tuple $(\alpha,\alpha',\eta,\eta',\beta,\beta')$ to be
the ortho-end of a group in the family.

As usual, we work one plane at a time.
That is, we look at the six-tuple $(\alpha,\alpha',n,n',b,b')$.
We have
\begin{prop} $(\alpha,\alpha',n,n',b,b')$ lies in the NSDC-$\theta_A$  planar family
of $(a,a',n,n',b,b')$ if and only if
$$ {\frac{(a+a')-(\alpha+\alpha')-i{\frac{(a-a')(|a-a'|)}{2}}\cdot \tan \theta_A}
{i(\alpha-\alpha')}}\in \IR.$$
\end{prop}

\begin{proof} We need to show that this condition is
equivalent to $a,a',\alpha, $ and $ \alpha'$ lying on the
pull back circle $\theta_A$.  Recall that the pull-back angle, the
pull-back distance, and the center of the pull-back circle  are
related  by
\begin{equation}\label{eq:t}t_{\theta_A} =
{\frac{|a-a'|}{2}}\cdot \tan \theta_A
\end{equation}
\begin{equation} \label{eq:center} c_{\theta_A}= c_{t_{\theta_A}}=
{ {\frac{a +a'}{2}}} + i{\frac{(a -a')(|a-a'|)}{4}}\cdot \tan
\theta_A = {\frac{a + a'}{2}}+i{\frac{(a-a')}{2}}t_{\theta_A}
\end{equation}

Equations (\ref{eq:t}) and (\ref{eq:center}) must also hold for
$\alpha, \alpha', \theta_{A'}$ and $t_{\theta_{A'}}$. The four
points lie on the same circle precisely when
$c_{\theta_A}=c_{\theta_{A'}}$. This happens if and only if there
is a real $t_{\theta_{A'}}$ with
$${ {\frac{a +a'}{2}}} + i{\frac{(a -a')(|a-a'|)}{4}}\cdot \tan
\theta_A = {\frac{\alpha
+\alpha'}{2}}+i{\frac{(\alpha-\alpha')}{2}}t_{\theta_{A'}}$$

\end{proof}


\begin{thebibliography}{99}
\bibitem{Fench} W. Fenchel, {\em Elementary Geometry in Hyperbolic Space},
de Gruyter Studies in Mathematics, $11$, Berlin-New York,(1989).
\bibitem{GiNSDC} J. Gilman, {\em A Discreteness condition for subgroups of
$PSL(2, \IC)$}, Proceedings of the Bers Colloquium, Contemporary Math Series,
$211$, AMS, (1997) 261-267.
\end{thebibliography}
\end{document}